\input amstex
\magnification=\magstep0
\documentstyle{amsppt}
\pagewidth{6.0in}
\input amstex
\topmatter
\title
Differential operators, shifted parts, \\
and hook lengths\\
\endtitle
\affil
Tewodros Amdeberhan\\ 
Massachusetts Institute of Technology, \\
Department of Mathematics, Cambridge, MA;\\
and \\
Tulane University, New Orleans, LA\\
\tt{tewodros\@math.mit.edu}
\endaffil
\abstract
We discuss Sekiguchi-type differential operators, their eigenvalues, and a generalization of Andrews-Goulden-Jackson formula.
These will be applied to extract explicit formulae involving shifted partitions and hook lengths. 
\endabstract
\endtopmatter
\def\({\left(}
\def\){\right)}

\document
\noindent
\bf 1. Differential operators. \rm
\smallskip
\noindent
The standard Jack symmetric polynomials $P_{\lambda}(y_1,\dots,y_n;\alpha)$ (see Macdoland, Stanley [5, 10]) 
as well as their shifted counter-parts (replace $\theta=1/\alpha$; see Okounkov-Olshanski [7] and references therein) have been studied. 
The former appear as eigenfunctions of the Sekiguchi differential operators
$$\align
&D(u;\theta)=a_{\delta}(y)^{-1}\det\left(y_i^{n-j}\left(y_i\frac{\partial}{\partial y_i}+(n-j)\theta+u\right)\right\}_{i,j=1}^n,\\
&D(u;\theta)P_{\lambda}(y;\theta)=\left(\prod_i(\lambda_i+(n-i)\theta+u)\right)P_{\lambda}(y;\theta),\tag1
\endalign$$
where $\delta:=(n-1,n-2,\dots,1,0)$ and $\lambda$ are partitions, $a_{\delta}=\prod_{1\leq i<j\leq n}(y_i-y_j)$ is the Vandermonde determinant 
and $u$ is a free parameter.
\bigskip
\noindent
Under a general result, S. Sahi proves [8, Theorem 5.2] the existence of a 
unique polynomial $P^{*}_{\mu}(y;\theta)$, now known as shifted 
Jack polynomials, satisfying
a certain \it vanishing condition. \rm In the special case $\theta=1$, 
Okounkov and Olshanski [6,7] relate $P^{*}_{\mu}(y;1)$ 
to the Schur functions $s_{\lambda}$ and call them \it shifted Schur 
polynomials. \rm
\bigskip
\noindent
One motivation for this paper is the a result due to Andrews, Goulden and Jackson [1, Thm 2.1] stating that
\smallskip
\noindent
\bf Theorem 1.1. \it For $n, m\in\Bbb{P}$ and summing over all partitions 
$l(\lambda)\leq n$, we have
$$\align
&\sum_{\lambda}s_{\lambda}(y_1,\dots,y_n)s_{\lambda}(w_1,\dots,w_m)\prod_{i=1}^n(x-\lambda_i-n+i)\\
&=\prod_{j=1}^n\prod_{k=1}^m(1-y_jw_k)^{-1}\times
[t_1,\dots,t_n](1+t_1+\cdots+t_n)^x\prod_{k=1}^m\left(1-\sum_{j=1}^n\frac{t_jy_jw_k}{1-y_jw_k}\right).
\endalign$$ \rm
\bf Remark 1.2: \rm To put things in perspective with the above operator view point, let us change $x$ to $-u$ and multiply through by $(-1)^n$ to find
$$\align
&G:=\sum_{\lambda}s_{\lambda}(y_1,\dots,y_n)s_{\lambda}(w_1,\dots,w_m)\prod_{i=1}^n(\lambda_i+(n-i)+u)\\
&=(-1)^n\prod_{j=1}^n\prod_{k=1}^m(1-y_jw_k)^{-1}\times
[t_1,\dots,t_n](1+t_1+\cdots+t_n)^{-u}\prod_{k=1}^m\left(1-\sum_{j=1}^n\frac{t_jy_jw_k}{1-y_jw_k}\right).\endalign$$ \rm
Observe that the differential operators in $D(u;\theta)$ do \it not \rm commute with the Vandermonde determinants, so the action may be
regraded as one-sided. On the other hand, a closer look at  Theorem 1.1 (rather its reformulation Remark 1.2) 
shows its basic underpinning being, in a sense, the opposite end of the bracket for Sakiguchi when $\theta=1$.
\bigskip
\noindent
\bf Definition 1.3. \rm A \it left-Sakiguchi operator \rm is the composed map $D^{\prime}(u;1):=a_{\delta}^{-1}L(u;1)\circ\psi_{\delta}$, where
$$L(u;1)=\prod_{i=1}^n(y_i\frac{\partial}{\partial y_i}+u),\qquad\text{and}\qquad \psi_{\delta}:F\rightarrow a_{\delta}F.$$ 
\bf Proposition 1.4. \it 
\smallskip
\noindent
(a) The Schur polynomials are eigenfunctions to $D^{\prime}(u;1)$; meaning that
$$D^{\prime}(u;1)s_{\lambda}=\left(\prod_i(\lambda_i+(n-i)+u)\right)s_{\lambda}.$$
In particular, the action of $D^{\prime}(u;1)$ on symmetric polynomials is diagonalizable having distinct eigenvalues.
\smallskip
\noindent
(b) The operators $\{D^{\prime}(u;1): u\in\Bbb{C}\}$ form a commutative algebra, say $\Cal{D}^{\prime}(n;1)$. Moreover,
if $\Lambda^{1}(n)$ is the algebra of polynomials $g(z_1,\dots,z_n)$ that are symmetric in the variables $z_i-i$ then
by part (a) we have the Harish-Chandra isomorphism
$$\Cal{D}^{\prime}(n;1)\rightarrow \Lambda^{1}(n)$$
mapping an operator $D$ to a polynomial $d(\lambda)$ such that $D(s_{\lambda})=d(\lambda)s_{\lambda}$.
\smallskip
\noindent
\bf Proof: \rm Since $a_{\delta}s_{\lambda}=a_{\delta+\lambda}$ it suffices to check that
$$a_{\delta+\lambda}\prod_{i=1}^n(\lambda_i+(n-i)+u)=\sum_{\sigma\in S_n}\text{sgn}(\sigma)\cdot 
\prod_{i=1}^n(\lambda_i+(n-i)+u)y_{\sigma(i)}^{\lambda_i+n-i}
=L(u;1)a_{\delta+\lambda}.$$
This however is straightforward. $\square$
\bigskip
\noindent
Now extend the left-Sakiguchi operation to 
$D^{\prime}(u;\theta):=a_{\theta\delta}^{-1}L(u;1)\circ \psi_{\theta\delta}$ acting
on symmetric functions. Denote $s^{\prime}_{\theta\delta+\lambda}=a_{\theta\delta}^{-1}a_{\theta\delta+\lambda}$ where $\theta\delta$ is
a simple \it homothety \rm map $(n-i)\rightarrow (n-i)\theta$. It is clear that
$$D^{\prime}(u;\theta)s^{\prime}_{\theta\delta+\lambda}=\left(\prod_i(\lambda_i+(n-1)\theta+u)\right)s^{\prime}_{\theta\delta+\lambda}.\tag2$$ \rm
Thus, both the Sakiguchi operators (1) and the full left-Sakiguchi (2) share the same eigenvalues.
One concludes $D(u;\theta)$ and $D^{\prime}(u;\theta)$ are similar transformations. Therefore, a parallel analysis regarding 
$D^{\prime}(u;1)$ can be carried out as in [4,6,7,8] but we defer such an undertaking to
the interested reader. Instead our main focus will be on the operator $L(u;1)$ and the eigenvalues
$\phi_{\lambda}(u):=\prod_i(\lambda_i+n-i+u)$. Keep in mind that
$$L(u;1)a_{\delta+\lambda}=\phi_{\lambda}(u)a_{\delta+\lambda}.$$
The paper is organized as follows.

\bigskip
\noindent
\bf 2. Shifted parts. \rm
\smallskip
\noindent
A conjecture on Guoniu Han [3, Conj. 3.1] asserts that for all positive integers $k$, the expression
$$\frac1{n!}\sum_{\lambda\vdash n}f_{\lambda}^2\sum_{v\in\lambda}h_{v}^{2k}$$
is a polynomial function of $n$. Richard P. Stanley actually settles this conjecture in
some generalized form [9, Theorem 4.3]. In the course of the proof he also show several intermediate results. 
In this section, we provide an alternative proof for one of those statements [9, Lemma 3.1].
\bigskip
\noindent
Suppose $\lambda=(\lambda_1,\dots,\lambda_n)$
is a partition of $n$, adding as many $0$'s at the tail as needed. Throughout this paper, $f_{\lambda}$ denotes
the number of standard Young tableaux (SYT) of shape $\lambda$ and $H_{\lambda}=\prod_{v\in\lambda}h_v$ will be
the product of the hook lengths $h_{v}$ where $v$ ranges over all squares in the Young diagram of $\lambda$. The content
of the cell $v$ is denoted $c_v$. Define
$$A_{\lambda}(u)=\frac{\phi_{\lambda}(u)}{H_{\lambda}}.$$
In the sequel, the following well-established facts (see Frame, Robinson and Thrall [2]) shall be untilized repeatedly
$$f_{\lambda}=\frac{n!}{H_{\lambda}}, \qquad \frac1{H_{\lambda}}=\det\left(\frac1{(\lambda_i-i+j)!}\right)_{i,j=1}^n.\tag 3$$
\bf Lemma 2.1 [9, Lemma 3.1] \it Let $p_1=y_1+\cdots+y_n$ and $e_i$ the elementary symmetric functions.
Then
$$\sum_{i=0}^n\binom{u+i-1}{i}p_1^ie_{n-i}=\sum_{\lambda\vdash n}A_{\lambda}(u)s_{\lambda}.$$
\bf Proof: \rm From Cauchy's $\sum_{\lambda}s_{\lambda}(w)s_{\lambda}(y)=\prod_{i,j=1}^n(1-w_iy_j)^{-1}$ 
we get $\sum_{\lambda\vdash n}f_{\lambda}s_{\lambda}(y)=p_1^n$.
We already know that $a_{\lambda+\delta}\phi_{\lambda}(u)=L(u;1)a_{\lambda+\delta}$. Therefore
$$\align
\sum_{\lambda\vdash n}f_{\lambda}s_{\lambda}(y)\phi_{\lambda}(u)
&=a_{\delta}^{-1}\sum_{\lambda\vdash n}f_{\lambda}a_{\lambda+\delta}\phi_{\lambda}(u)\\
&=a_{\delta}^{-1}L(u;1)\sum_{\lambda\vdash n}f_{\lambda}a_{\lambda+\delta}\\
&=a_{\delta}^{-1}L(u;1)a_{\delta}\sum_{\lambda\vdash n}f_{\lambda}s_{\lambda}\\
&=a_{\delta}^{-1}L(u;1)a_{\delta}\cdot p_1^n.\endalign$$
A direct computation using Leibnitz' Rule for multi-derivatives of products yields
$$L(u;1)a_{\delta}p_1^n=a_{\delta}\cdot n!\sum_{j=0}^n\binom{u+j-1}{u-1}p_1^je_{n-j}.$$
Consequently, we have
$$\sum_{\lambda\vdash n}f_{\lambda}s_{\lambda}\phi_{\lambda}(u)=a_{\delta}^{-1}L(u;1)a_{\delta}p_1^n
=n!\sum_{i=0}^n\binom{u+i-1}{i}p_1^ie_{n-i}.$$ 
By the hook length formula $f_{\lambda}=n!H_{\lambda}^{-1}$ the proof is complete. $\square$
\bigskip
\noindent
\bf Corollary 2.2. \it Let $u\in\Bbb{N}, \lambda\vdash n, \mu:=(n^u,\lambda_1^{\prime},\lambda_2^{\prime},\dots,\lambda_n^{\prime})\vdash (u+1)n$
and $\lambda^{\prime}=(\lambda_1^{\prime},\cdots,\lambda_n^{\prime})$ the shape conjugate to $\lambda$. Then the quantity $A_{\lambda}(u)$
\smallskip
\noindent
(a) enumerates the product of the hook lengths of the top row in $\mu$ divided by $H_{\lambda^{\prime}}$ or $H_{\lambda}$; symbolically
$A_{\lambda}(u)=\frac{H_{\mu}}{H_{\lambda}^2}$.
\smallskip
\noindent
(b) is an integer via the following combinatorial identity involving weighted skew SYT
$$A_{\lambda}(u)=\sum_{i=0}^n\binom{u+n-i-1}{n-i}f^{\lambda/1^i}.$$ \rm
\bf Proof: \rm The symmetric polynomials $p_1^{n-i}e_i$ can be expanded in terms of the Schur basis $\{s_{\lambda}\}_{\lambda\vdash n}$
as follows: $p_1^{n-i}e_i=\sum_{\lambda\vdash n}f^{\lambda/1^i}s_{\lambda}$. Now apply the above lemma to obtain
$$\sum_{\lambda\vdash n}A_{\lambda}(u)s_{\lambda}=\sum_{\lambda\vdash n}s_{\lambda}\sum_{i=0}^n\binom{u+n-i-1}{n-i}f^{\lambda/1^i}.$$ 
Comparison of coefficients yields the identity of part (b). The assertion of part (a) is evident. $\square$
\bigskip
\noindent
\bf Remark 2.3. \rm Although $H_{\lambda}^{-1}\prod_i(\lambda_i+(n-i)+u)$ is an integer, the quantity
$$\frac1{H_{\lambda}}\prod_{i=1}^n(\lambda_i+(n-i)+u_i)$$
is \it not \rm generally an integer even if all $u_i$'s are. Take for instance, $n=2, \lambda=(2,0), u_1=1$ and $u_2=2$.
Thus $H_{\lambda}^{-1}(3+u_1)u_2=\frac12(5)(1)\notin\Bbb{Z}$.
\bigskip
\noindent
\bf 3. Determinants. \rm
\smallskip
\noindent
In this section, we compile a few determinantal evaluations some of which are well-known and others are residues from
the previous sections. The first statement is valid for arbirary parameters $z_i$'s and the proof is tailor-made for
Dodgson's method of \it condensation \rm (for a charming proof, see [11]).
\bigskip
\noindent
\bf Proposition 3.1. \it Given the variables $z=(z_1,\dots,z_n)$, there holds
$$\det\left(\binom{z_i+n-i}{n-j}\right)_{i,j=1}^n
=\prod_{1\leq i<j\leq n}\frac{z_i-z_j+j-i}{j-i}.$$ \rm
\bf Proof: \rm Introduce two integral parameters $b, c$ to generalize the claim further as
$$\det\left(\binom{z_{i+b}+n-(i+b)}{n-(j+c)}\right)_{i,j=1}^n
=\prod_{1\leq i<j\leq n}\frac{z_{i+b}-z_{c+j}+(j+c)-(i+b)}{(j+c)-(i+b)}.$$
The rest is an automatic application of Dodgson. $\square$
\bigskip
\noindent
\bf Corollary 3.2. \it Let $\lambda\vdash n$ with length $l=l(\lambda)$ and 
$\omega_i(u)=\lambda_i+(n-i)+u$.
$$\det\left(\frac1{(\lambda_i-i+j)!}\right)_{i,j=1}^n
=\frac{\prod_{1\leq i<j\leq n}(\lambda_i-\lambda_j+j-i)}
{\prod_{i=1}^n(\lambda_i+n-i)!}=\frac1{H_{\lambda}}.\tag a$$
$$\det\left(\frac{i}{(\lambda_i-i+j)!}\right)_{i,j=1}^n=\frac{n!}{H_{\lambda}}=f_{\lambda}\in\Bbb{Z}.\tag b$$
$$\det\left(\frac{\omega_i(u)}{\omega_i(j-n)!}\right)_{i,j=1}^n
=\prod_{i=1}^n\frac{\omega_i(u)}{\omega_i(0)!}\cdot
\prod_{i<j}\omega_i(u)-\omega_j(u)=A_{\lambda}(u)\in\Bbb{Z};\tag c$$
with the convention that $\frac1{y!}:=0$ whenever $y<0$.
$$\det\left(\binom{\lambda_i+l-i}{l-i}\right)_{i,j=1}^l=\frac1{H_{\lambda}}\prod_{i=1}^l\frac{(\lambda_i+l-i)!}{(l-i)!}=
\prod_{v\in\lambda}\frac{l+c_u}{h_u}\in\Bbb{Z}.\tag d$$ \rm
\bf Proof: \rm In Prop. 3.1, divide by $(z_i+n-i)!(n-j)!^{-1}$ then replace $\lambda_i$ for $z_i$. This together with (3) prove part(a).
The remaining parts follow directly from Prop. 3.1 and Cor. 2.2. $\square$

\bigskip
\noindent
\bf 4. Andrews-Goulden-Jackson formula. \rm
\smallskip
\noindent
We begin with an extraction of certain coefficients from Theorem 1.1 (or equivalently Remark 1.2). Notation: given a
function $F(x,y,\dots)$, the coefficient of the monomial $x^ay^b\cdots$ is desiganted by $[x^ay^b\cdots]F$. Also $s(n,k)$ and $c(n,k)$ stand
for the signed and unsigned stirling numbers of the first kind, respectively. The falling factorials $(x)_k=x(x-1)\cdots(x-k+1)$ generate
$$(x)_k=\sum_{i=0}^ks(k,i)x^i.$$
\bf Lemma 4.1. \it The coefficient of the falling factorial $(x)_{n-k}$ in Remark 1.2 is
$$\prod_{i=1}^n\prod_{j=1}^m(1-y_iw_j)^{-1}
\sum\Sb j_1<\cdots<j_k\\i_1\neq\cdots\neq i_k\endSb\prod_{r=1}^k\frac{y_{j_1}w_{i_r}}{1-y_{j_1}w_{i_r}}.$$ 
Moreover, if we extract the coefficient of $[w_1,\dots,w_n]$ we get
$$(n)_k\cdot(y_1+\cdots+y_n)^{n-k}\cdot
\sum\Sb j_1<\cdots<j_k\endSb\prod_{r=1}^ky_{j_1}=(n)_k\cdot(y_1+\cdots+y_n)^{n-k}\cdot e_k(y_1,\cdots,y_n).$$ \rm
\bf Corollary 4.2. \it Let $m=n$ and fix $\beta$. Then the following 
coefficients are polynomials in $n$, of degree $2\beta$. 
$$\frac1{n!}[y_1\cdots y_nw_1\cdots w_nx^{n-\beta}]G=\frac1{n!}\sum_{\lambda\vdash n}f_{\lambda}^2e_{\beta}(\{\lambda_i+n-i:1\leq i\leq n\})=
\sum_{\alpha=n-\beta}^nc(\alpha,n-\beta)\binom{n}{\alpha}.$$ \rm
\bf Examples 4.3: \rm We list the first few polynomials $\frac1{n!}\sum_{\lambda\vdash n}f_{\lambda}^2e_{\beta}$.
\bigskip
\noindent
a) for $\beta=0$ we have $\binom{n}0$.
\bigskip
\noindent
b) for $\beta=1$ we have $\binom{n+1}2$.
\bigskip
\noindent
c) for $\beta=2$ we have $-\binom{n}2-\binom{n+1}3+3\binom{n+2}4$.
\bigskip
\noindent
d) for $\beta=3$ we have 
$\binom{n}3-5\binom{n+1}4-10\binom{n+2}5+15\binom{n+3}6$.
\bigskip
\noindent
e) for $\beta=4$ we have $2\binom{n}4+19\binom{n+1}5-20\binom{n+2}6-105\binom{n+3}7+105\binom{n+4}8$.

\bigskip
\noindent
\bf 5. Main Results: Towards a generalization. \rm
\bigskip
\noindent
Let $y=(y_1,\dots,y_n), w=(w_1,\dots,w_n)$. Recall the eigenvalues and eigenfunctions of $D^{\prime}(u;1)$:
$$\phi_{\lambda}(u)=\prod_i(\lambda_i+(n-i)+u), \qquad
L(u;1)=\prod_i(y_i\frac{\partial}{\partial y_i}+u).$$
For brevity, write $L(u)$ instead of $L(u;1)$. Next we proceed to compute
$$\sum_{\lambda}s_{\lambda}(w)s_{\lambda}(y)\prod_{j=1}^N\phi_{\lambda}(u_j)
=\sum_{\lambda}s_{\lambda}(w)s_{\lambda}(y)\prod_{j=1}^N\prod_{i=1}^n(\lambda_i+(n-i)+u_j).\tag{STAN-N}$$
Note once again that
$$a_{\delta+\lambda}\prod_{j=1}^N\phi_{\lambda}(u_j)=\sum_{\sigma\in S_n}\text{sgn}(\sigma)\cdot 
\prod_{j=1}^N\phi_{\lambda}(u_j)y_{\sigma(i)}^{\lambda_i+n-i}
=\prod_{j=1}^NL(u_j)a_{\delta+\lambda}.$$
This in turn implies
$$\align
\sum_{\lambda}s_{\lambda}(w)s_{\lambda}(y)\prod_{j=1}^N\phi_{\lambda}(u_j)
&=a_{\delta}^{-1}(y)\sum_{\lambda}s_{\lambda}(w)a_{\delta+\lambda}(y)
\prod_{j=1}^N\phi_{\lambda}(u_j)\\
&=a_{\delta}^{-1}(y)\left(\prod_{j=1}^NL(u_j)\right)\sum_{\lambda}s_{\lambda}(w)
a_{\delta+\lambda}(y)\\
&=a_{\delta}^{-1}(y)\left(\prod_{j=1}^NL(u_j)\right)a_{\delta}(y)\sum_{\lambda}
s_{\lambda}(w)s_{\lambda}(y)\\
&=a_{\delta}^{-1}(y)\left(\prod_{j=1}^NL(u_j)\right)a_{\delta}(y)\prod_{i,k=1}^n(1-y_iw_k)^{-1},\endalign$$
where the last equality uses Cauchy's summation formula. 
\bigskip
\noindent
\bf Case $N=1$: \rm When $m=n$ this is exactly Remark 1.2 that will be re-formulated as:
$$\align
&a_{\delta}^{-1}(y)L(u_1)a_{\delta}(y)\prod_{i,k=1}^n(1-y_iw_k)^{-1}\\
&=(-1)^n\prod_{j=1}^n\prod_{k=1}^n(1-y_jw_k)^{-1}\times
[t_1,\dots,t_n](1+t_1+\cdots+t_n)^{-u_1}\prod_{k=1}^n\left(1-\sum_{j=1}^n\frac{t_jy_jw_k}{1-y_jw_k}\right).\endalign$$
If we denote the right-hand side by $\bold{RHS}(u_1)$ then
$$\align
\sum_{\lambda}s_{\lambda}(w)s_{\lambda}(y)\prod_{j=1}^N\phi_{\lambda}(u_j)
&=a_{\delta}^{-1}(y)L(u_N)\cdots L(u_2)a_{\delta}\left(a_{\delta}^{-1}L(u_1)a_{\delta}\prod_{i,k=1}^n(1-y_iw_k)^{-1}\right)\\
&=a_{\delta}^{-1}(y)L(u_N)\cdots L(u_2)a_{\delta}\text{$\bold{RHS}$}(u_1).\endalign$$
\bf Case $N=2$: \rm We opt to proceed where we left-off in Lemma 4.1, i.e. we work on the coefficients of $[w_1\cdots w_n](u_1)_{n-k}$.
The advantage is two-fold: the notational amount is reduced and it is irrelevant to maintain the extra family of parameters $w_i$.
$$(n)_kL(u_2)a_{\delta}(y_1+\cdots+y_n)^{n-k}e_k(y_1,\dots,y_n).\tag4$$
Since $(n-k)![s^k][z^n]e^{zp_1(y)}\prod_{i=1}^n(1+y_isz)=p_1^{n-k}e_k(y)$, where $p_1(y)=y_1+\cdots+y_n$,
one may consider instead
$$(n)_kL(u_2)a_{\delta}F(y_1)\cdots F(y_n), \qquad\text{with}\qquad  F(\alpha)=e^{z\alpha}(1+\alpha sz).$$
At this point we invoke [1, Theorem 4.3] of Andrews, Goulden and Jackson asserting that
$$\align
L(u_2)a_{\delta}F(y_1)\cdots F(y_n)
&=(-1)^na_{\delta}F(y_1)\cdots F(y_n)[t_1\cdots t_n](1+t_1+\cdots+t_n)^{-u_2}\\
&\exp\left\{\left(\alpha\frac{\partial}{\partial\alpha}\log 
F(\alpha)\right)*\log\left(1-\sum_{j=1}^n\frac{t_jy_j\alpha}{1-y_j\alpha}\right)\right\}.
\endalign$$
Here $\sum_{i\geq0}a_iz^i *\sum_{i\geq0}b_iz^i=\sum_{i\geq0}a_ib_i$. 
Note distributivity $(a+b)*c=a*c+b*c$. For our choice of the function 
$F(\alpha)$ it is then clear that
$$\align
\alpha\partial_{\alpha}\log F(\alpha)&=z\alpha+\frac{\alpha sz}{1+\alpha sz},\\
\exp\left\{(z\alpha)*
\log\left(1-\sum_{j=1}^n\frac{t_jy_j\alpha}{1-y_j\alpha}\right)\right\}&=\exp\left\{-z\sum_{j=1}^nt_jy_j\right\},\\
\exp\left\{\left(\frac{\alpha sz}{1+\alpha sz}\right)*
\log\left(1-\sum_{j=1}^n\frac{t_jy_j\alpha}{1-y_j\alpha}\right)\right\}&=\left(1+\sum_{j=1}^n\frac{t_jy_jsz}{1+y_jsz}\right)^{-1}.\endalign$$
Therefore
$$\align
a_{\delta}^{-1}(n)_kL(u_2)a_{\delta}e^{zp_1(y)}\prod_{i=1}^n(1+y_isz)
&=(-1)^n(n)_ke^{zp_1(y)}\prod_{i=1}^n(1+y_isz)[t_1\cdots t_n](1+t_1+\cdots+t_n)^{-u_2}\\
&e^{-zp_1(\bold{ty)}}\sum_{j=0}^n(-1)^jj!(sz)^je_j(\bold{ty});\tag5\endalign$$
where $p_1(\bold{ty})=t_1y_1+\cdots+t_ny_n$ and 
$e_j(\bold{ty})=e_j(t_1y_1,\cdots,t_ny_n)$. 
\bigskip
\noindent
The key is now to extract the coefficient of $(n-k)![y_1\cdots y_n][s^k][z^n]$ from both sides of equation (5), and
the matter rests on what we find on the right-hand side. The procedure has four stages.
\bigskip
\noindent
\bf Coefficients $[s^k]$: \rm
$$[s^k]\prod_{i=1}^n(1+y_isz)\cdot\sum_{j=0}^n(-1)^jj!(sz)^je_j(\bold{ty})
=z^k\sum_{\alpha=0}^k(-1)^{\alpha}\alpha!e_{\alpha}(\bold{ty})e_{k-\alpha}(y).$$
\bf Coefficients $[z^n]$: \rm Since there is already a gain of $z^k$, actually we must read off $[z^{n-k}]$. So,
$$[z^{n-k}]e^{zp_1(y)}e^{-zp_1(\bold{ty})}=\sum_{\beta=0}^{n-k}\frac{(-1)^{\beta}}{\beta!(n-k-\beta)!}p_1(y)^{n-k-\beta}p_1(\bold{ty})^{\beta}.$$
\bf Coefficients $[t_1\cdots t_n]$: \rm
$$[t_1\cdots t_n](1+t_1+\cdots 
t_n)^{-u_2}p_1(\bold{ty})^{\beta}e_{\alpha}(\bold{ty})
=\frac{(\alpha+\beta)!}{\alpha!}e_{\alpha+\beta}(y)(-u_2)_{n-\alpha-\beta}.$$
\bf Coefficients $[y_1\cdots y_n]$: \rm
$$[y_1\cdots y_n]p_1(y)^{n-k-\beta}e_{k-\alpha}(y)e_{\alpha+\beta}(y)=
(n-k-\beta)!\binom{n}{\alpha+\beta}\binom{n-\alpha-\beta}{k-\alpha}.$$
It is now time to combine all these into
$$\align
(n-k)!(n)_k &\sum_{\alpha=0}^k\sum_{\beta=0}^{n-k}
\frac{(-1)^{\alpha+\beta}\alpha!(\alpha+\beta)!(n-k-\beta)!}
{\beta!(n-k-\beta)!\alpha!}
\binom{n}{\alpha+\beta}\binom{n-\alpha-\beta}{k-\alpha}(-u_2)_{n-\alpha-\beta}\\
&=n!\sum_{\alpha=0}^k\sum_{\beta=0}^{n-k}
\frac{(-1)^{\alpha+\beta}(\alpha+\beta)!}{\beta!}
\binom{n}{\alpha+\beta}\binom{n-\alpha-\beta}{k-\alpha}(-u_2)_{n-\alpha-\beta}\\
&=n!\sum_{\gamma=0}^n(-1)^{\gamma}\gamma!\binom{n}{\gamma}(-u_2)_{n-\gamma}
\sum_{\alpha=0}^k\frac1{(\gamma-\alpha)!}\binom{n-\gamma}{k-\alpha}.
\endalign$$
\bf Lemma 5.1 (STAN $N=2$). \it Fix $\alpha, \beta\in\Bbb{P}$. Then, the coefficient of $u_1^{n-\alpha}u_2^{n-\beta}$ is given by
$$\align
&\frac1{n!}\sum_{\lambda\vdash n}f_{\lambda}^2e_{\alpha}(\{\lambda_i+n-i:1\leq i\leq n\})e_{\beta}(\{\lambda_j+n-j:1\leq j\leq n\})\\
&=\sum_{k=n-\alpha}^n\sum_{m=n-\beta}^nc(k,n-\alpha)\cdot c(m,n-\beta)\cdot\binom{n}m\sum_{j\geq 0}
j!\cdot\binom{n-m}j\binom{m}{n-k-j}.\endalign$$
Here, the index $j$ runs through its natural limits, i.e. in the sense that if $Y > X \geq 0$ then $\binom{X}{Y}=0$. \rm
\bigskip
\noindent
\bf Proof: \rm Derivation seen above. $\square$
\bigskip
\noindent
Let $K=(k_{ij})$ be an infinite lower triangular matrix of non-negative integral entries, and define the row-sum and differences by
$$\bar{r}_{i}=k_{i1}+\cdots+k_{ii}, \qquad \text{and} \qquad
\Delta \bar{r}_{i-1}=\bar{r}_{i-1}-\bar{r}_i+k_{ii}.$$
\bf Theorem 5.2. \it  For any integer $N\geq 1$, we have 
$$\align
\frac1{n!}\sum_{\lambda\vdash n}s_{\lambda}f_{\lambda}\prod_{j=1}^N\phi_{\lambda}(u_j)
=\sum_{\epsilon=1}^N&\sum_{k_{\epsilon\epsilon}=0}^n(n-k_{\epsilon\epsilon})!\binom{u_{\epsilon}+n-k_{\epsilon\epsilon}-1}{n-k_{\epsilon\epsilon}}\times\\
&\sum_{\alpha=2}^N\sum \Sb\Delta\bar{r}_{\alpha-1}\geq 0\endSb \binom{n-k_{NN}}{\Delta\bar{r}_N-k_{NN}}\Delta\bar{r}_N
\binom{k_{\alpha\alpha}}{\Delta\bar{r}_{\alpha-1}}\Delta\bar{r}_{\alpha-1}p_1^{n-\bar{r}_N}\prod_{i=1}^Ne_{k_{N,i}}.\endalign$$ \rm
\bf Proof: \rm A successive application of the argument that produces Lemma 5.1. $\square$
\bigskip
\noindent
Another one of Stanley's results [9, Theorem 3.3] proves polynomiality of the
expression
$$\Psi_n(W)=\frac1{n!}\sum_{\lambda\vdash n}f_{\lambda}^2W(\{\lambda_i+n-i:1\leq i\leq n\};\{c_v: v\in\lambda\})\tag6$$
where $W(x,y)$ is a given power series over $\Bbb{Q}$ of bounded degree
that is symmetric in the variables $x=(x_1,x_2,\dots)$ and $y=(y_1,y_2,\dots)$
separately. In fact, the main kernel of our motivation for this paper emanates from studying $\Psi_n(W)$.
\bigskip
\noindent
Our main result, the corollary below, only deals with the 
shifted parts $\{\lambda_i+n-i:1\leq i\leq n\}$ that appear in equation (6) and not the 
contents $c_v$. However the ouput is indeed an explicit evaluation, the
method of proof is different and the conclusion of \it polynomiality \rm is preserved. 
\bigskip
\noindent
\bf Corollary 5.3 (STAN-N). \it  Given the $N$-tuple diagonal $(k_{11},\dots,k_{NN})$ of $K$, we have
$$\align
&\frac1{n!}\sum_{\lambda\vdash n}f^2_{\lambda}
\prod_{j=1}^Ne_{k_{jj}}(\{\lambda_i+n-i: 1\leq i\leq n\})\\
&=\sum_{\gamma=1}^N\sum_{j_{\gamma \gamma}=n-k_{\gamma\gamma}}^n
c(j_{\gamma\gamma},n-k_{\gamma\gamma})
\sum_{\alpha=2}^N\sum\Sb \Delta\bar{r}_{\alpha-1}\geq 0\endSb 
(k_{\alpha\alpha})_{\Delta\bar{r}_{\alpha-1}}\binom{n}{k_{N1},\dots,k_{NN},n-\bar{r}_N}.\endalign$$
If $W(x)$ is a power series over $\Bbb{Q}$ of bounded degree and symmetric in $x=(x_1,x_2,\dots)$, then
$$\frac1{n!}\sum_{\lambda}f_{\lambda}^2W(\{\lambda_i+n-i: 1\leq i\leq n\})$$ 
is a polynomial function of $n$. \rm
\bigskip
\noindent
\bf Proof: \rm Apply Theorem 5.2 and the generating function $(x)_k=\sum_{i=0}^ks(k,i)x^i$ for stirling numbers, to get the formula.
In particular polynomiality follows. $\square$

\bigskip
\noindent
\bf Acknowledgments. \rm The author is grateful to Rosena Du for introducing him to Guoniu Han's recent work 
on hook lengths. The author also thanks Richard P. Stanley for many
inspiring discussions and his valuable ideas.

\enddocument

\noindent
\bf Question for Stanley: \rm
\bigskip
\noindent
Thank you for showing me that
$$p_1^{n-i}e_i=\sum_{\lambda\vdash n}f^{\lambda/1^i}s_{\lambda}.$$
\it Here is a question now: is the following true? If not, what is the right-hand side? \rm
$$p_1^{n-k_1-k_2-\dots -k_r}e_{k_1}e_{k_2}\cdots e_{k_r}=
\sum_{\lambda\vdash n}f^{\lambda/\mu}s_{\lambda}, \qquad \mu=(1^{k_1},1^{k_2},\dots,1^{k_r}).$$
The understanding with $\lambda/\mu$ is deleting $k_1$ boxes from the first column of $\lambda$,
then deleting $k_2$ boxes from the second column of $\lambda$, and so on.


$$\frac{(n+l)!}{l!} \prod_{i=1}^l\frac{(l-i)!}{(\lambda_i+l-i)!}$$

\noindent


\bigskip
\noindent
Next, introduce the operator $L(\xi)=\prod_{i=1}^n(\xi-y_i\partial_{y_i})$. Keep in mind that $L(\xi_s)$ and $L(\xi_u)$ do not commmute. 
Then, in principle, the next \it reduction \rm statement generalizes the Theorem 2.1.
\bigskip
\noindent
\bf Lemma 4: \it Let $\xi_1,\dots,\xi_N$ be indeterminantes. We have
$$\sum_{\theta}s_{\theta}(y_1,\dots,y_n)s_{\theta}(w_1,\dots,w_m)\prod_{\alpha=1}^N\prod_{i=1}^n(\xi_{\alpha}-\theta_i-n+i)
=a_{\delta}^{-1}L(\xi_N)\cdots L(\xi_2)\left\{a_{\delta}\times\text{\bf RHS($\xi_1$)}\right\};$$ \it
where by \bf{RHS($\xi_1$)} \it we mean the right-hand side of Theorem 2.1 treated as a function of $\xi_1$. \rm
\bigskip
\noindent
\bf Remark 5: \rm The leading term on the right-side of the Lemma 4 is
$$\text{\bf{RHS}($\xi_1$)$\times$\bf{RHS}($\xi_2$)$\times\cdots\times$\bf{RHS}($\xi_N$)}.$$ \rm
\bf Corollary 6: \it Denote the falling factorial by $(\xi)_n=\xi(\xi-1)\cdots(\xi-n+1)$. In the special case $m=0$, it holds
$$\sum_{\theta}s_{\theta}(y_1,\dots,y_n)\prod_{\alpha=1}^N\prod_{i=1}^n(\xi_{\alpha}-\theta_i-n+i)=(\xi_1)_n(\xi_2)_n\cdots(\xi_N)_n.$$
Less trivially,
$$a_{\delta}^{-1}L(\xi_N)\cdots L(\xi_1)a_{\delta}=(\xi_N)_n\cdots(\xi_1)_n.$$ \rm
\bigskip
\noindent
Consider the case $N=2$. The coefficient of $(x)_k$ in $G$ is
$$\prod_{i=1}^n\prod_{j=1}^n(1-y_iw_j)^{-1}\cdot
\sum\Sb j_1<\cdots< j_k\\i_1\neq\cdots\neq i_k\endSb\prod_{l=1}^k\frac{y_{j_l}w_{i_l}}{1-y_{j_l}w_{i_l}}
=\tilde{A}\cdot B.$$
Assume $m\leq n$. Let $A=a_{\delta}\tilde{A}$, $\partial_{y_i}\rightarrow\partial_i$, $L_{i}(\xi)=\xi-y_i\partial_i$ and 
$-y_i\partial_{y_i}B\rightarrow\partial_iB$ then we compute
$$L(\xi)(A\cdot B)=\sum\Sb 1\leq i_1<\cdots<i_m\leq n\\1\leq i_{m+1}<\cdots<i_n\leq n\endSb
(L_{i_1}\cdots L_{i_m}A)\cdot(\partial_{i_{m+1}}\cdots\partial_{i_n}B).$$

\bigskip
\noindent
To model the operations on $A$, it suffices to study $\{i_1,\dots,i_m\}=[m]$. Also let $M=2n-m$ and define
$$\omega_k(j)=\cases \frac{-y_jw_k}{1-y_jw_k},\qquad{} \qquad1\leq k\leq n; \\
\frac{y_j}{y_{k+m-n}-y_j},\qquad n+1\leq k\leq M. \endcases$$
It follows that the coefficient of $(\xi)_r$ in $a_{\delta}^{-1}L_1\cdots L_mA$ takes the form
$$\prod_{i=1}^n\prod_{j=1}^n(1-y_iw_j)^{-1}\cdot
\sum\Sb j_1<\cdots< j_r\\i_1\neq\cdots\neq i_r\endSb\prod_{l=1}^r\omega_{i_l}(j_l);$$
where the summation runs through $j_q\in[m]$ and $i_q\in[M]$.
\bigskip
\noindent
On the other hand, $\partial_{m+1}\dots\partial_nB$ vanishes unless  $\{m+1,\dots,n\}\subset\{j_1,\cdots,j_k\}$ (hence
$n-m\leq k$), in which case the value becomes (here $i_s\in[n]$, $y_{i_l}\neq \tilde{y_l}=y_{l+n-k}$)
$$\sum\Sb j_1<\cdots<j_{k+m-n}\\i_1\neq\cdots\neq i_k\endSb\prod_{l=1}^{k+m-n}
\frac{y_{j_l}w_{i_l}}{1-y_{j_l}w_{i_l}}\prod_{l=k+m-n+1}^k\frac{-\tilde{y_l}w_{i_l}}{(1-\tilde{y_l}w_{i_l})^2}.$$


An eye on extendability of Theorem 2.1. Throughout we work in the field $\Bbb{F}_2$.
\bigskip
\noindent
\bf Case-study 1: \it Suppose $\bold{x}=(x_1,\dots,x_m)$ and let 
$$f_n(\bold{x}):=\left(1+\sum_{i=1}^mx_i+x_1\sum_{i=2}^mx_i^2\right)^{2^n-1}\in\Bbb{F}_2[\bold{x}].$$
Then, the number of odd coefficients is
$$N(f_n(\bold{x}))=m\cdot(m+1)^n-(m-1)\cdot m^n.$$ \rm
\bf Proof: \rm Working in $\Bbb{F}_2[\bold{x}]$ and applying Fermat's $(a_1+\cdots+a_k)^p\equiv(a_1^p+\cdots+a_k^p)$, it is evident that 
$$f_n(\bold{x})=\prod_{j=0}^{n-1}(1+\sum_{i=0}^mx_i^{2^j}+\sum_{i=2}^mx_1^{2^j}x_i^{2^{j+1}})
=\prod_{j=0}^{n-1}(1+\sum_{i=2}^mx_i^{2^j}+\cdots+\sum_{i=1}^m\bar{x}_i^{2^j});$$
where $\bar{x}_1:=x_1\cdot1^2$ and $\bar{x}_i:=x_1\cdot x_{i}^2$ for $0\leq i\leq m$. There's a certain benefit in the
following rearrangment which involves rotating every other row
$$\align
f_n(\bold{x})=&(1+x_2+\cdots+x_m+\bar{x}_1+\bar{x}_2+\cdots+\bar{x}_m)\cdot\\
&(\bar{x}_m^2+\cdots+\bar{x}_1^2+\bar{x}_1^2+x_m^2+\cdots+x_2^2+1^2)\cdot\\
&(1^{4}+x_2^{4}+\cdots+x_m^{4}+\bar{x}_1^{4}+\bar{x}_2^{4}+\cdots+\bar{x}_m^{4})\cdots\qquad\text{\& etc}.\endalign$$
We intend to show that $b_n:=f_n(\bold{x})$ satisfies a three-term recurrence, and it suffices to explain the idea
with the three consecutive products  $b_{n-2}, b_{n-1}, b_n$ (the reader can work on $b_1, b_2, b_3$). 
Note that in principle each row has $2m$ summands.
\bigskip
\noindent
Consider $b_{n-1}$: the list of duplicate \it directed edges, \rm resulting from multiplying the two rows, can be coded by 
$(\bar{x}_i\rightarrow x_j^2)=(\bar{x}_j\rightarrow x_i^2)$
provided $i\neq j$. This happens on either the west or east-coast $(2\times m)$-block (for definiteness, say east). 
That means each of the $m$ NE-nodes has out-degree 
$m-1$ (the same in-degree at each SE-nodes). On the other hand, in the second row there are $m$ SW-nodes each receiving the full 
$2m$ in-degree (or $b_{n-2}$). Hence 
$$b_{n-1}=mb_{n-2}+m(b_{n-2}-(m-1)b_{n-3}).\tag1$$
Focus on $b_n$: the same argument as above is applicable between the second and third rows, and on the west-coastal $(2\times m)$-block. 
However, there are additional interactions between the $1^{st}$ and $3^{rd}$ rows, via the $2^{nd}$. In fact, the duplication
occurs between the NE $m$-block (first row) and SW $m$-block (third row) where the NE $m$-block (second row) acts as a catalyst.
The process is to construct paths from the first row through the single nodes we left out in $b_{n-1}$ into the second row
and proceed as in $b_{n-1}$ but into the third row. The code is 
$(\bar{x}_i\rightarrow x_i^2\rightarrow x_j^2)=(\bar{x}_j\rightarrow1\rightarrow x_i^4)$, for $i\neq j$.
This list must be deleted, so
$$b_n=mb_{n-1}+m(b_{n-1}-(m-1)b_{n-2}-(m-1)b_{n-3}).\tag2$$
Combining (1) and (2) leads to
$$b_n=(2m+1)b_{n-1}-m(m+1)b_{n-2}.$$
We arrive at the assertion after solving this recurrence with $b_0=1, b_1=2m$. $\square$

\bigskip
\noindent

\bigskip
\noindent
\bf Examples: \rm
\bigskip
\noindent
a) If $m=2$ then $f_2(x,y)=(1+x+y+xy^2)^{2^n-1}\in\Bbb{F}_2[x,y]$ and $N(f_2(x,y))=2\cdot3^n-2^n$.
The same formula holds for $g_n=(1+x+y^2+xy)^{2^n-1}$.
\bigskip
\noindent
b) Denote forward shift $E(a_n)=a_{n+1}$. Below, we add a few more illustrations in $\Bbb{F}_2[\bold{x}]$.
\bigskip
\noindent
i) When $m=3$, $f_n=(1+x+y+z+xy^2+xz^2)^{2^n-1}$ and $a_n=N(f_n)$ satisfies $E^2-7E+12=(E-4)(E-3)=0$. So, $a_n=3\cdot4^n-2\cdot3^n$. 
\bigskip
\noindent
ii) If $g_n=(1+x+y+z+xy^2+yz^2)^{2^n-1}$ then $a_n=N(g_n)$ satisfies $E^2-7E+10=(E-5)(E-2)=0$. So, $a_n=\frac13(4\cdot5^n-2^n)$.
\bigskip
\noindent
iii) Let $g_n=(1+x+y+z+xy^2)^{2^n-1}$, then $a_n=N(g_n)$ satisfies $E^2-6E+7=0$. Here the eigenvalues are not even rational; i.e.
if $\lambda=3+\sqrt{2}$ and $b=1+\sqrt{2}$ and write $\bar{\mu}$ for quadratic conjugation in $\Bbb{Z}[\sqrt{2}]$ then
$$a_n=\frac{b\lambda^n+\bar{b}\bar{\lambda}^n}2.$$
iv) Let $g_n=(1+x+y+z+xy+xz^2)^{2^n-1}$. Now $a_n=N(g_n)$ satisfies $E^2-7E+8=0$. Again the eigenvalues are not even rational; i.e.
if $\lambda=\frac{7+\sqrt{17}}2$ and $b=17+5\sqrt{17}$ and write $\bar{\mu}$ for quadratic conjugation in $\Bbb{Q}[\sqrt{17}]$ then
$$a_n=\frac{b\lambda^n+\bar{b}\bar{\lambda}^n}{34}.$$
v) Let $g_n=(1+x+y^2+xy^3)^{2^n-1}$, $\lambda=2+\sqrt{3}$ and conjugates in the quadratic extension $\Bbb{Z}[\sqrt{2}]$. We have
$$N(g_n)=\frac{\lambda^{n+1}+\bar{\lambda}^{n+1}-2^n}3.$$
vi) (\bf Unsolved\rm). If $g_n=(1+x+y+x^2y^3)^{2^n-1}$ and $g^{\prime}_n=(1+x^2+y+xy^2)^{2^n-1}$, then $N(g_n)=N(g^{\prime}_n)$.
\bigskip
\noindent
c) \it Periodic solutions. \rm 
\smallskip
\noindent
i) Let $g_n=(1+x+y+x^2y^2)^{2^n-1}$ then $N(g_n)$ satisfies $E^4-5E^3+6E^2-2E-4=0$ whose eigenvalues are complex. 
Put $\lambda=\frac{3+\sqrt{17}}2$ and $b=17^2+73\sqrt{17}$ with conjugates in $\Bbb{Q}[\sqrt{17}]$. Thus
$$N(g_n)=\frac{b\lambda^n+\bar{b}\bar{\lambda}^n}{442}-\frac{2}{\sqrt{13}}2^{n/2}\cos\left(\frac{\pi n}4-\arctan(3/2)\right).$$
ii) Let $g_n=(1+x^2+y^2+xy^3)^{2^n-1}$ then there is a generating function for $N(g_n)$ given by
$$\sum_{n\geq 0}N(g_n)\xi^n=\frac{1-2\xi+4\xi^2}{1-6\xi+12\xi^2-12\xi^3}.$$
iii) Let $g_n=(1+x+y+y^2)^{2^n-1}$ then 
$$\sum_{n\geq 0}N(g_n)\xi^n=\frac{1+2\xi}{1-2\xi-4\xi^2}.$$
In fact $N(g_n)=2^nF(n+2)$ with $F(n)$ fibonacci. The following is also immediate.
\smallskip
\noindent
\bf Corollary: \it 
\smallskip
\noindent
1) Given $k\in\Bbb{P}$, break up the binary digits of 
$k$ into maximal strings of consecutive $1$'s and let the lengths of these
strings be the set $\bold{k}=\{k_1,k_2,\dots\}$. Then we obtain the 
averaging value
$$\frac1{2^n}\sum_{k=0}^{2^n-1}\prod_{\bold{k}}\frac{2^{k_i+2}+(-1)^{k_i+1}}3
=F(n+2).$$ 
2) The probability of not landing two consecutive heads in a fair toss of $n$ coins is equal to that of finding
a head in the $(2^n-1)$-triangle formed by $(1+y+y^2)^m$ mod $2$, which is $F(n+2)/2^n$. \rm (Head=$1$, Tail=$0$.)
\bigskip
\noindent
\it 3) The generating function $\Lambda(\xi)=\sum_{m=0}^{\infty}N((1+y+y^2)^m)\xi^m$ satisfies $\Lambda(\xi)=(1+2\xi)\Lambda(\xi^2)$. \rm
\bigskip
\noindent
iv) Let $g_n=(1+x+y+y^3)^{2^n-1}$ then 
$$\sum_{n\geq 0}N(g_n)\xi^n=\frac{1+\xi-2\xi^3}{1-3\xi-2\xi^2+2\xi^3+4\xi^4}.$$
Further, the generating function $\Lambda(\xi)=\sum_{m=0}^{\infty}N((1+y+y^3)^m)\xi^m$ satisfies the equation
$$\align
\Lambda(\xi)&=a(\xi)\Lambda(\xi^2)+b(\xi)\Lambda(\xi^4)+c(\xi)\Lambda(\xi^8)+d(\xi)\Lambda(\xi^{16});\qquad\text{where}\\
a(\xi)&=1+2\xi^8, \qquad{} \qquad{} \qquad c(\xi)=-\xi^2+4\xi^3+2\xi^6-6\xi^{10}-2\xi^{12}+\xi^{14}\\
b(\xi)&=3\xi+\xi^2-2\xi^8, \qquad \qquad d(\xi)=-2\xi^6+2\xi^{12}-9\xi^{14}-2\xi^{15}+7\xi^{16}.\endalign$$
v) Let $g_n=(1+x+y+y^4)^{2^n-1}$ then  
$$\sum_{n\geq 0}N(g_n)\xi^n=\frac{1+\xi+4\xi^2+2\xi^3-4\xi^4}{1-3\xi-2\xi^3-8\xi^4+8\xi^5}.$$
vi) If $g_n=(1+x+y+z+xy^2+xz^2+yz^2)^{2^n-1}$ then $N(g_n)$ is generated by
$$\sum_{n\geq 0}N(g_n)\xi^n=\frac{(1-\xi)^2}{1-9\xi+23\xi^2-19\xi^3}.$$
vii) If $g_n=(1+x+y+z+xy^2+yx^2)^{2^n-1}$ then $N(g_n)$ satisfies 
$$\sum_{n\geq 0}N(g_n)\xi^n=\frac{1-\xi+2\xi^2-4\xi^3}{1-7\xi+12\xi^2-12\xi^3+8\xi^4}.$$
The degree of complexity increases 
if the polynomials we seek are even small deviants from the above special form $f_n$; this is more so when some
measure of symmetry is lacking.
\bigskip
\noindent
d) \it Symmetric polynomials. \rm Consider the Vandermonde-type 
polynomials, in $\Bbb{F}_2[\bold{x}]$, given by
$$V(k,n):=\prod_{1\leq i<j\leq k}(x_i-x_j)^{2^n-1} \qquad\text{and}\qquad
V^{\prime}(k,n):=\left(1+\prod_{1\leq i<j\leq k}(x_i-x_j)\right)^{2^n-1}.$$
We find that
\bigskip
\noindent
i) $N(V(2,n))=2^n$ and $N(V^{\prime}(2,n))=3^n$.
\bigskip
\noindent
ii) If $\lambda=\frac{5+\sqrt{33}}2, b=11+3\sqrt{33}$ with conjugation in $\Bbb{Q}[\sqrt{33}]$ then 
$$N(V(3,n))=6\cdot4^{n-1}, \qquad N(V^{\prime}(3,n))=\frac{b\lambda^n+\bar{b}\bar{\lambda}^n}{22}.$$
Moreover, if $M=\pmatrix 1&2q\\\frac3q&4\endpmatrix$ or $M=\pmatrix 1&q\\\frac6q&4\endpmatrix$ then $N(V^{\prime}(3,n))=$ the first entry of $M^{n+1}$.
\bigskip
\noindent
iii) $N(V(4,n))=5\cdot8^n-8\cdot2^n$.
\bigskip
\noindent
e) For more, check Example 3.1 of our paper.

\bigskip
\noindent
\bf Formula 5.2. \rm

\bigskip
\noindent
First some nomenclature. Suppose $\lambda=(\lambda_1,\dots,\lambda_k)$ is a partition of $n\in\Bbb{Z}^{+}$ so that $\sum_{i=1}^k\lambda_i=n$ and 
$0<\lambda_i$'s are weakly decreasing. Write $\lambda\vdash n$. The set of all partitions of $n$ is designated $\Cal{P}(n)$ while the total
number of partitions of $n$ will be $p(n)$. 
\smallskip
\noindent
Each partition $\lambda$ is represented by its Ferrer's diagram and the corresponding Young Tableaux of shape $\lambda$. For
each box $v\in\lambda$ in the diagram associate the hook length $h_v(\lambda)$ (or simply $h_v$). Let $f^{\lambda}$ stand for 
the number of standard Young Tableaux of shape $\lambda$. Given any complex number $\alpha$, introduce an \it $\alpha$-norm 
\rm 
$\Vert\lambda\Vert^{\alpha}=\lambda_1^{\alpha}+\cdots+\lambda_k^{\alpha}$. 
Then,
\smallskip
\noindent
I) the log-average over $\Cal{P}(n)$ becomes
$$\frac1{p(n)}\sum_{\lambda\vdash n}\sum_{v\in\lambda}\log\left(f^{\lambda}h_v^n\right)=n\log(n!).$$
True, this is immediate from \it hook formula. \rm
\bigskip
\noindent
II) the following overall sum on $\Cal{P}(n)$ holds
$$\sum_{\lambda\vdash n}\sum_{v\in\lambda}h_v^{\alpha-1}=
\sum_{\lambda\vdash n}\Vert\lambda\Vert^{\alpha},\qquad
\sum_{\lambda\vdash n}\sum_{v\in\lambda}q^{h_v}=\sum_{\lambda\vdash n}\sum_i\lambda_iq^{\lambda_i}.$$
III) there is a generating function for (II) (both II and III were found unaware of Han's results, too bad)
$$\frac1{(z;z)_{\infty}}\sum_{k\geq1}\frac{z^kk^{\alpha}}{(z;z)_k}=
\frac{\sum_{k\geq1}\sigma_{\alpha}(k)z^k}{(z;z)_{\infty}}=
\frac{\zeta(z)\zeta(z-\alpha)}{(z;z)_{\infty}}.$$ 
Here $(z;z)_{\omega}=\prod_{j=1}^{\omega}(1-z^j)$;
$\sigma_{\alpha}(k)=\sum_{d\vert k}d^{\alpha}$ the 
$\alpha^{th}$-divisor function; $\zeta(z)$ Riemann's zeta. 
\bigskip
\noindent
IV) if $l(\lambda)$ is the length of $\lambda$, it is easy to find
$$\sum_{n\geq 0}\left(\sum_{\lambda\vdash 
n}\prod_{i=1}^{l(\lambda)}\lambda_i^{\alpha}\right)z^n=\prod_{j\geq1}
\frac1{1-j^{\alpha}z^j}.$$
\bf Corollary 1: \it Suppose $F(z)$ is an entire function. Then, we have
$$\align \sum_{\lambda\vdash n}\sum_{v\in\lambda}\frac{F(h_v)}{h_v}&=\sum_{\lambda\vdash n}\sum_iF(\lambda_i),\\
\sum_{\lambda\vdash n}\sum_{v\in\lambda}F(h_v)&=\sum_{\lambda\vdash n}\sum_i\lambda_iF(\lambda_i).\endalign$$ \rm
\bf Corollary 2: \it It holds that 
$$\align
\prod_{\lambda\vdash n}\prod_{h\in\lambda}h&=\prod_{\lambda\vdash n}\prod_i\lambda_i^{\lambda_i},\\
\sum_{\lambda\vdash n}\sum_{h\in\lambda}h^{h-1}&=\sum_{\lambda\vdash n}\sum_i\lambda_i^{\lambda_i},\\
\sum_{\lambda\vdash n}\sum_{h\in\lambda}\left(1+\frac1{h^{\alpha}}\right)^{\beta}
&=\sum_{\lambda\vdash n}\sum_i\lambda_i\left(1+\frac1{\lambda_i^{\alpha}}\right)^{\beta}.\endalign$$
\bf Corollary 3: \it Suppose $\Cal{X}_n$ is the character table of the symmetric group $S_n$, $m_i(\lambda)=$ multiplicities 
of $i$ in $\lambda$. Then we have
$$\prod_{\lambda\vdash n}\prod_{h\in\lambda}\root h\of{h}=\vert\det(\Cal{X}_n)\vert=\prod_{\lambda\vdash n}\prod_i\lambda_i
=\prod_{\lambda\vdash n}\prod_{i}m_i(\lambda)!.$$ \rm
\bf Proof: \rm All corollaries are immediate from (II). $\square$
\bigskip
\noindent
In [1, page 11], Conjecture 5.2 states that
$$\sum_{\lambda\in\Cal{P}}x^{\vert\lambda\vert}
\prod_{h\in\Cal{H}(\lambda)}\rho(z;h)=e^{x+zx^2/2},\tag A$$
with the weight function $\rho(z;n)$ is defined by
$$\rho(z;n)=\frac{\sum_{k=0}^{\lfloor{n/2}\rfloor}\binom{n}{2k}z^k}
{n\sum_{k=0}^{\lfloor{(n-1)/2}\rfloor}\binom{n}{2k+1}z^k}$$
and $h$ is the hook-length.
\bigskip
\noindent
It is perhaps more facile to have $\rho$ re-expressed as (we replaced $z\rightarrow z^2$):
$$\rho(z^2;n)=\frac{z}n\left(\frac{(1+z)^n+(1-z)^n}{(1+z)^n-(1-z)^n}\right).$$
Then equation (A) takes the following equivalent format
$$\align
e^{x+z^2x^2/2}
&=\sum_{\lambda\in\Cal{P}}x^{\vert\lambda\vert}\prod_{h\in\Cal{H}(\lambda)}\frac{z}h
\left(\frac{(1+z)^h+(1-z)^h}{(1+z)^h-(1-z)^h}\right)\\
&=\sum_{\lambda\in\Cal{P}}(zx)^{\vert\lambda\vert}\prod_{h\in\Cal{H}(\lambda)}\frac1h
\left(\frac{(1+z)^h+(1-z)^h}{(1+z)^h-(1-z)^h}\right)\\
&=\sum_{\lambda\in\Cal{P}}\frac{(zx)^{\vert\lambda\vert}}{\vert\lambda\vert!}f_{\lambda}\prod_{h\in\Cal{H}(\lambda)}
\frac{(1+z)^h+(1-z)^h}{(1+z)^h-(1-z)^h}\\
&=\sum_{\lambda\in\Cal{P}}\frac{(zx)^{\vert\lambda\vert}}{\vert\lambda\vert!}f_{\lambda}\prod_{h\in\Cal{H}(\lambda)}
\coth\left(\frac{h}2\log\left(\frac{1+z}{1-z}\right)\right);
\endalign$$
where the third equality uses the \it hook formula \rm and $\coth$ is the 
hyperbolic cotangent. Still, if we substitute $z\rightarrow \tanh\omega$ then we obtain
$$\sum_{\lambda\in\Cal{P}}\frac{(x\tanh\omega)^{\vert\lambda\vert}}{\vert\lambda\vert!}f_{\lambda}\prod_{h\in\Cal{H}(\lambda)}
\coth(h\omega)=e^{x+(x\tanh\omega)^2/2};\qquad or,$$
$$\sum_{\lambda\in\Cal{P}}x^{\vert\lambda\vert}\prod_{h\in\Cal{H}(\lambda)}
\frac1h\frac{\tanh\omega}{\tanh(h\omega)}=e^{x+(x\tanh\omega)^2/2}.$$
\bf Conjecture: \it In lieu of the last expression, we project that similar statements would hold for
a class of weight functions $\Psi$, with appropriate boundary conditions (at $\omega=0$ and $\omega\rightarrow\infty$), so that
$$\sum_{\lambda\in\Cal{P}}x^{\vert\lambda\vert}\prod_{h\in\Cal{H}(\lambda)}
\frac1h\frac{\Psi(\omega)}{\Psi(h\omega)}$$
can be given in terms of elementary (or classical) functions. \rm Test: $\Psi(\omega)=\tan^{-1}\omega$.
\bigskip
\noindent
\bf Towards $\sum\prod\frac1{h^2}\sum h^{2m}$ \rm
\bigskip
\noindent
Given a positive integer $m$, consider all its partitions $\alpha=(\alpha_1,\dots,\alpha_k)\vdash m$. Begin with the expansion 
$$\left(\sum_hh^2\right)^m=\sum_{\alpha\vdash m}\prod_{(u_1,\dots,u_k)}h_{u_i}^{2\alpha_i}.\tag *$$
The inner-most product runs through all \it ordered \rm $k$-tuples with $u_i\neq u_j$ whenever $i\neq j$.
\bigskip
\noindent
\bf Some observations: \rm
\smallskip
\noindent \it
a) There exist a polynomial $P\in\Bbb{Z}[x]$ of degree $2m$ and a constant $\gamma\in\Bbb{Q}$ depending only on $m$ such that
$$\sum_{\lambda\vdash n}\prod_{v\in\lambda}\frac1{h_v^2}\left(\sum_{u\in\lambda}h_u^2\right)^m=\frac{\gamma_m}{n!}P_{2m}(n).$$
b) Consider each term in (*). Then, there exist a polynomial $Q\in\Bbb{Z}[x]$ of degree $m$ and a constant $\epsilon\in\Bbb{Q}
$ depending only on $k,m$ so that
$$\sum_{\lambda\vdash n}\prod_{v\in\lambda}\frac1{h_v^2}\sum_{(u_1,\dots,u_k)}\prod_{i}h_{u_i}^{2\alpha_i}=\frac{\epsilon_{k,m}}{(n-k)!}Q_m(n).$$ 
c)  If $(n)_j$ denotes the falling factorial $n(n-1)\cdots(n-j+1)$, then RHS of (a)-(b) take the form
$$\align
\frac1{n!}\sum_{\lambda\vdash n}f_{\lambda}^2\left(\sum_{u\in\lambda}h_u^2\right)^m&=\gamma_m\sum_{j=1}^{2m}\sigma_{j,m}(n)_j,\\
\frac1{n!}\sum_{\lambda\vdash n}f_{\lambda}^2\sum_{(u_1,\dots,u_k)}\prod_{i}h_{u_i}^{2\alpha_i}&=\epsilon_{k,m}\sum_{j=k}^{k+m}\tilde{\sigma}_{j,m}(n)_j;\endalign$$ 
and the coefficient sequences $\sigma, \tilde{\sigma}$ are \bf positive \it integers with \bf no internal \it zeros and \bf unimodal\rm. 

\bigskip
\noindent
\bf Example: \rm
\bigskip
\noindent
The case $m=3$. All relevant quantities involved are computed below.
$$\align
\frac1{n!}\sum_{\lambda\vdash n}f_{\lambda}^2\left(\sum_{u\in\lambda}h_u^2\right)^3&=
\frac{24(n)_1+1476(n)_2+5720(n)_3+4746(n)_4+1170(n)_5+81(n)_6}{24},\\
\frac1{n!}\sum_{\lambda\vdash n}f_{\lambda}^2\sum_{u\in\lambda}h_u^6&=
\frac{24(n)_1+756(n)_2+2240(n)_3+1050(n)_4}{24},\\
\frac1{n!}\sum_{\lambda\vdash n}f_{\lambda}^2\sum_{(u_1,u_2)}h_{u_1}^4h_{u_2}^2&=
\frac{60(n)_2+254(n)_3+235(n)_4+60(n)_5}6,\\
\frac1{n!}\sum_{\lambda\vdash n}f_{\lambda}^2\sum_{(u_1,u_2,u_3)}h_{u_1}^2h_{u_2}^2h_{u_3}^2&=
\frac{144(n)_3+292(n)_4+150(n)_5+27(n)_6}8.
\endalign$$
\it Caveat: \rm  The number $\bold{552}$ has been corrected from Han's papers. We count ordered tuples while he counts \it unordered \rm ones
(the conversion is trivial). The $1^{st}$ and the $3^{rd}$ calculations are new. 
\bigskip
\noindent
\bf Miscellaneous: \rm
\bigskip
\noindent \it
a) Let $f(\lambda):=\#\{v\in\lambda: \text{arm$(v)$-leg$(v)$ $=0,1$}\}$. Then
$$\align        
\sum_{\lambda\vdash n}f(\lambda)^{\beta}&=\sum_{\lambda\vdash n}l(\lambda)^{\beta},\\
\sum_{\lambda\vdash n}\binom{f(\lambda)}{\beta}&=\sum_{\lambda\vdash n}\binom{l(\lambda)}{\beta},\\
\sum_{\lambda\vdash n}\binom{\beta+f(\lambda)}{\beta}^{-1}&=\sum_{\lambda\vdash n}\binom{\beta+l(\lambda)}{\beta}^{-1}.\endalign$$
b) If $p_{n,k}=$ the number of partitions of $n$ into exactly $k$ parts, then
$$\sum_{\lambda\vdash n}q^{f(\lambda)}=\sum_{\lambda\vdash n}q^{l(\lambda)}=\sum_{k=1}^np_{n,k}q^k;\qquad
p_{n,k}=[x^n]\frac{x^k}{(x;x)_k}.$$
c) Define the sequence $b_1=b_2=1$ and $b_n=\frac{b_{n-2}+b_{n-1}}{n-2}$, for $n\geq 3$. Then, there is an integer $\gamma_k$ such that
$$e^{-3/2}\sum_{n\geq k}\binom{n}kb_n=\gamma_k.$$
The first few values: $\{\gamma_0,\gamma_1,\dots\}=\{1,3,3^2,3^3+2,3^4+12,\dots\}$, The sequence $b_n$ arose as a coefficient in
$$[w^2x^n]e^{x+(x\cosh\omega)^2/2}=\frac12b_{n-1}.$$

\enddocument